\documentclass[a4paper,11pt]{article}
\pdfoutput=1 % if your are submitting a pdflatex (i.e. if you have
\usepackage{jheppub-ac} % for details on the use of the package, please

\AtBeginDocument{%
  \addtolength\abovedisplayskip{-0.2\baselineskip}%
  \addtolength\belowdisplayskip{-0.2\baselineskip}%
}

\usepackage[T1]{fontenc} % if needed

\usepackage{comment}
\usepackage{upgreek}
\usepackage{bm}
\usepackage[normalem]{ulem}
\usepackage{mathrsfs}
\usepackage{amsfonts}

\newcommand{\nb}[1]{{\color{blue} NB: #1}}
\newcommand{\ac}[1]{\noindent\textcolor{red}
{{\rm [\![}\mbox{\sc{AC}$\blacktriangleright\!\!\blacktriangleright$}: {#1}{\rm ]\!]}}}

\renewcommand{\ac}[1]{}
\newcommand{\ak}[1]{}
\renewcommand{\nb}[1]{}

\newcommand{\calA}{\mathcal{A}}
\newcommand{\calL}{\mathcal{L}}
\newcommand{\scrX}{\mathscr{X}}
\newcommand{\scrY}{\mathscr{Y}}
\newcommand{\frakL}{\mathfrak{L}}
\newcommand{\bmK}{\bm{K}}
\newcommand{\jj}{\mathrm{i}}
\newcommand{\e}{\bm{e}}

\newcommand{\hh}{Y}
\DeclareSymbolFont{EUR}{U}{eur}{m}{n}
\SetSymbolFont{EUR}{bold}{U}{eur}{b}{n}
\DeclareSymbolFontAlphabet{\eur}{EUR}
\DeclareSymbolFont{EUB}{U}{eur}{b}{n}
\SetSymbolFont{EUB}{bold}{U}{eur}{b}{n}
\DeclareSymbolFontAlphabet{\eub}{EUB}
\newcommand{\GG}{\eub{G}}
\renewcommand{\GG}{\bm\Upomega}
\newcommand{\p}{\partial}
\DeclareSymbolFont{AMSb}{U}{msb}{m}{n}
\DeclareSymbolFontAlphabet{\mathbb}{AMSb}
\newcommand{\R}{\mathbb{R}}\newcommand{\C}{\mathbb{C}}
\newcommand{\N}{\mathbb{N}}

\newcommand{\abs}[1]{\vert #1 \vert}

\newcommand{\norm}[1]{\Vert #1 \Vert}
\renewcommand{\Re}{\mathop{\rm{R\hskip -1pt e}}\nolimits}

\title{Stable bi-frequency spinor modes as Dark Matter candidates}
\author[a,*]{Andrew Comech,\note[*]{Corresponding author.}}
\author[a]{Niranjana Kulkarni,}
\author[b]{Nabile Boussa\"{i}d,}
\author[c,d]{Jes\'us Cuevas-Maraver}

\affiliation[a]{Mathematics Department, Texas A\&M University, College Station, TX 77843-3368, USA}
\affiliation[b]{
UTINAM (UMR 6213), Universit\'e Marie et Louis Pasteur, \'equipe $\varPhi$th,
  \\F-25000 Besan\c{c}on, France}

\affiliation[c]{Grupo de F\'{\i}sica No Lineal, Departamento de F\'{\i}sica Aplicada I,
  Universidad de Sevilla,\\Escuela Polit\'{e}cnica Superior,
  C/ Virgen de Africa, 7, 41011-Sevilla, Spain
}

\affiliation[d]{Instituto de Matem\'{a}ticas de la Universidad de Sevilla (IMUS),\\
  Edificio Celestino Mutis. Avda. Reina Mercedes s/n, 41012-Sevilla, Spain
}

\emailAdd{comech@tamu.edu}
\emailAdd{nkulkarni@tamu.edu}
\emailAdd{nabile.boussaid@univ-fcomte.fr}
\emailAdd{jcuevas@us.es}

\abstract{
We show that
spinor systems with scalar self-interaction,
such as the Dirac--Klein--Gordon system
with Yukawa coupling
or the Soler model,
generically have
bi-frequency solitary wave solutions.
We develop the approach to stability properties of such waves
and use the radial reduction
to show that indeed the (linear) stability is available
for a wide range of parameters.
We show that only bi-frequency modes can be dynamically stable
and suggest that stable bi-frequency modes can serve as storages
of the Dark Matter.
The approach is based on linear stability results
of one-frequency solitary waves in (3+1)D Soler model,
which we obtain as a by-product.
}


\arxivnumber{2501.04027}

\usepackage{atbegshi}
\usepackage{fancyhdr}

\fancypagestyle{autorizzazione} {%
   \fancyfoot[C]{
 \footnotesize
 {\sc Open Access,} $\scriptstyle\copyright$ The Authors.
 \hfill
 \url{https://doi.org/10.1007/JHEP10(2025)082}
 \\
 \footnotesize
 Article funded by SCOAP$^3$.
 \hfill
}}

\begin{document}
\thispagestyle{autorizzazione}
\maketitle
\flushbottom

\tableofcontents

\setcounter{page}{1}
\section{Introduction}

Localized modes -- or quasiparticles --
are well-known in the classical field theories.
%%(that is, before the second quantization is applied).
These include polarons from condensed matter physics
%(see e.g. \cite{alexandrov2008polarons})
and skyrmions,
topological solitons in nonlinear sigma models.
%initially used to describe nucleons
%\cite{zahed1986skyrme},
%(see e.g. \cite{alexandrov2008polarons,everschor-sitte2018perspective}).
Polarons
%Both of these types of localized solutions
%enjoy stability properties
are related to
physical phenomena
such as
charge transport, surface reactivity,
colossal magnetoresistance, thermoelectricity,
photoemission, (multi)ferroism,
and high-temperature superconductivity
\cite{franchini2021polarons}.
%Polarons were initially used to describe nucleons;
%\cite{zahed1986skyrme},
%(see e.g.
%\cite{alexandrov2008polarons,everschor-sitte2018perspective}),
Magnetic skyrmions,
discovered in 2009
\cite{muuehlbauer2009skyrmion},
%are particle-like topological magnetic textures that
are now under consideration
as potential information carriers in spintronics
%\cite{finocchio2016magnetic,fert2017magnetic,li2023magnetic}.
\cite{li2023magnetic}.
%At the same time,
On the other hand,
localized modes of classical spinor fields
would always be treated with certain prejudice.
Indeed,
the Dirac sea hypothesis, the one which
prohibits electrons from descending into
negative energy states,
is based on the second quantization and
the Pauli exclusion principle, and it would seem to fail
for classical spinor fields,
supposedly rendering them unstable
and ready to plunge into the negative energy states.
In spite of this,
the nonlinear Dirac equation (NLD)
%% of the form
%% \[
%% \Big(\jj\gamma^\nu\frac{\p}{\p x^\nu}+g(\bar\psi\psi)\Big)\psi=0
%% \]
was considered by Ivanenko \cite{jetp.8.260}
%as a model where an electron-positron pair
%is created not by some heavier particle
%but by electron itself
and then
by Finkelstein and others and by Heisenberg
\cite{finkelstein1951nonlinear,finkelstein1956nonlinear,heisenberg1957quantum}
as a model of relativistic quantum matter.
The NLD appears
in the Nambu--Jona--Lasinio model in the hadron theory
\cite{nambu1961dynamical},
in the theory of Bose--Einstein condensates
\cite{PhysRevLett.104.073603},
and
%% as effective equations 
in photonics
\cite{smirnova2020nonlinear}.
Nonlinear spinor models
are discussed in the context
of Quantum Gravity, Cosmology,
Dark Matter, and Dark Energy \cite{wang2016dark}.
%In 
%\cite{dzhunushaliev2021mass,dzhunushaliev2024charged},
%nonlinear configurations
%of spinor field
%coupled to the Yang--Mills and electromagnetic fields
%are considered as models of elementary particles.
%with the attention to the existence of
%\emph{mass gap}; see also \cite{dzhunushaliev2019non}.
%We mention here, though, that
%in a rather general situation
%the smallest energy nonlinear modes
%are exactly on the border of stability and instability regions
%of parameters;
%these ``critical'' modes themselves
%are unstable \cite{comech2003purely}.

%Yet, systems of classical spinor fields
%have stable modes \cite{berkolaiko2012spectral,lakoba2018numerical},
%which makes such models physically viable.
%%\begin{comment}
To be physically viable,
a configuration of the fields needs to be stable;
there were numerous empirical attempts to address
stability of classical self-interacting spinor modes
as early as in the fifties.
%% The presumed instability
%% of spinor modes
%% was already questioned in
%% \cite{ranada1983classical}
%% on the grounds that the Dirac equation, which
%% is of first order in time, may fail our physical
%% intuition based on Newton law that the systems
%% tend into lower energy states.
It was suggested
%by Finkelstein et al. and then by Soler
\cite{finkelstein1951nonlinear,PhysRevD.1.2766}
that the smallest energy solitary wave
might be stable
(and then shown
%then Alvarez and Soler showed that
it was not \cite{PhysRevD.34.644};
as the matter of fact,
the linearization at
the minimal energy solitary wave
is characterized by
the collision of eigenvalues at zero
\cite{grillakis1987stability,boussaid2019nonlinear}
and consequently is unstable \cite{comech2003purely}).
Besides the numerical simulations
\cite{alvarez1981interaction,alvarez1983numerical}
which suggested stability in particular cases,
there were attempts to address
stability of spinor modes
based on energy or energy vs. charge considerations,
in the spirit of the energy approach
by Derrick \cite{derrick1964comments}
%(developed for the nonlinear wave equation)
and
the Grillakis--Shatah--Strauss theory \cite{grillakis1987stability};
we mention
\cite{%bogolubsky1979spinor,
%PhysRevD.34.644,
strauss1986stability}.
It was finally demonstrated
that spinor modes do possess
stability properties
for certain values of parameters,
on the examples of the (massive) Gross--Neveu model
and the (generalized) massive Thirring model
\cite{berkolaiko2012spectral,pelinovsky2014orbital,berkolaiko2015vakhitov,MR3462129,lakoba2018numerical}
and (2+1)D  Soler model
\cite{PhysRevLett.116.214101}.
%\cite{PhysRevD.12.3880}
\begin{comment}
As the matter of fact, stability
in nonlinear spinor systems
may be more robust than in more common scalar models:
%seem to demonstrate rather peculiar stability properties:
in particular, the blow-up behavior, which is well-established
in certain cases
for the nonlinear Schr\"odinger equation,
is presently unknown for the NLD
(under the assumption that the self-interaction
preserves $\mathbf{U}(1)$-invariance).
\end{comment}

Further studies of the NLD
\cite{boussaid2018spectral}
revealed a phenomenon
intrinsic to systems of
% spinor and scalar fields or
spinors with scalar self-interaction:
besides
%% usual, one-frequency
``Schr\"odinger-type'' modes
\begin{align}\label{of}
\psi_{\omega}(t,x)=\upvarphi(x)e^{-\jj\omega t}\in\C^4,
\qquad
x\in\R^3,
\quad
\omega\in\R,
\end{align}
which are known to exist in the NLD
since \cite{PhysRevD.1.2766},
such systems admit localized
bi-frequency modes of the form
\begin{equation}\label{bf-0}
\varPsi_{\pm\omega}(t,x)=
a\upvarphi(x)e^{-\jj\omega t}+b\upchi(x)e^{\jj\omega t},
\qquad
a,\,b\in\C,
\quad
\abs{a}^2-\abs{b}^2=1,
\end{equation}
with certain spatially localized $\upvarphi,\,\upchi$.
%(see \eqref{bf-1} below).
%% these modes exist
%% on the classical level, without introducing the second
%% quantization.
The phenomenon of bi-frequency modes
has been overlooked for years,
in spite of the discovery
\cite{galindo1977remarkable}
of $\mathbf{SU}(1,1)$ symmetry
in the Dirac--Klein--Gordon system (DKG)
and in the NLD:
\begin{align}\label{galindo}
\psi\mapsto
(a+b\mathsf{C})\psi
=
\left(a+\jj b
\gamma^2
\bmK\right)\psi,
\qquad
a,\,b\in\C,
\quad
\abs{a}^2-\abs{b}^2=1,
\end{align}
with $\gamma^2=
\Big[
\begin{matrix}0&\sigma_2\\[-4pt]-\sigma_2&0
\end{matrix}
\Big]
$
the corresponding Dirac matrix,
$\bmK$ the complex conjugation,
and $\mathsf{C}=\jj\gamma^2\bmK$ the charge conjugation operator;
one can see that the transformation \eqref{galindo}
yields bi-frequency modes
\eqref{bf-0} from \eqref{of}.
Most interestingly, though, is that
bi-frequency modes (see \eqref{bf-1} below)
are generically of more general form
than can be obtained via transformations \eqref{galindo}
(except in spatial dimensions $\le 2$ \cite{boussaid2018spectral});
\begin{comment}
or described
as the action of a one-parameter subgroup
of the symmetry group on a stationary state
(unlike traveling waves $u(x-ct)$ in the KdV
theory or standing waves
$e^{-\jj\omega t}\upvarphi(x)$
to the nonlinear Schr\"odinger equation),
\end{comment}
%%not fitting into the
their stability does not follow
from the
Grillakis--Shatah--Strauss stability theory
of standing waves
\cite{grillakis1990stability}
which is applicable to solutions of the form
$e^{\Omega t}\upvarphi$,
with $\Omega$ the Lie algebra
of the corresponding symmetry group
and $\upvarphi$ time-independent and localized in space.
% the study of linear stability of
The approach to stability of bi-frequency modes
has been absent.
%is to be developed from scratch.

Let us emphasize that
it is only bi-frequency modes that
can be \emph{dynamically}
(asymptotically) stable:
a bi-frequency mode \eqref{bf-0} with $\abs{b}\ll 1$,
considered as a small perturbation of \eqref{of},
cannot converge to a one-frequency mode,
since it is itself an exact solution.
We conclude that it is bi-frequency modes,
not one-frequency ones,
which may be of particular
interest for potential applications.
Dynamically stable bi-frequency modes
\eqref{bf-0}
can then provide models
for phenomena involving stable localized states
in the framework of spinor fields.

%% The coupling of the spinor field
%% to the (real-valued) Klein--Gordon field
%% can be interpreted as
%% interaction of fermions with Higgs bosons.
%% (this reminds of difficulties encountered by
%% ensembles of rotating charges in J.J. Thomson's plum--pudding model \cite{thomson1904structure}).
%Scalar-type interaction with the Klein--Gordon field
Yukawa-type interaction
in the DKG system
(the
$g\phi\bar\psi\psi$ term in the Lagrangian)
suggests that
bi-frequency modes can be considered
in relation to Dark Matter (DM)
theory (see e.g. \cite{arcadi2020dark}),
which is presently in search of suitable
candidates for DM particles:
%A well-motivated candidate for particle dark matter is
%with a Weakly Interacting Massive Particle (WIMP)
%being one of such candidates.
%(also in relation to the Dark Matter);
%An important possibility is that
stable neutral bi-frequency spinor modes
in the DKG system
%interacting only with the Higgs field
%-- the DKG system --
can model massive particles in the DM sector
interacting with the observed matter via the
``Higgs portal'',
as discussed in
\cite{cao2009dark,boehmer2010dark,bai2014lepton}.
Let us mention that
models of spinor-based DM are rather popular
\cite{bahamonde2018dynamical},
particularly so the ELKO spinors
%(\emph{Eigenspinoren des Ladungskonjugationsoperators})
\cite{da2011exotic,alves2015searching}.
%%\cite{alves2015searching}.
%(let us also mention $\mathbf{O}(3)$ spinors
%\cite{kitabayashi2020spinorial}).
We show below that classical bi-frequency modes
can be arbitrarily large
while retaining their stability properties,
which makes them possible storages of DM.

Bi-frequency modes,
interpreted as a particle-antiparticle superposition,
may also model other phenomena related to DM,
such as neutron--mirror neutron oscillations
$n\mbox{-}n'$
\cite{berezhiani2006neutron,berezhiani2009more,kamyshkov2022neutron,broussard2022experimental,dvali2024kaluza}
(with mirror neutron $n'$ considered
to be from the DM sector),
\mbox{neutron lifetime anomaly \cite{berezhiani2019neutron},}
\newpage
physics of neutron stars \cite{goldman2022neutron},
and sterile neutrino oscillations
\cite{boyarsky2016sterile,dasgupta2021sterile}.
They can also model
neutron--antineutron oscillations
$n\mbox{-}\bar n$
\cite{phillips2016neutron}.
%All these phenomena
%are presently under scrutiny in High Energy Physics
%as it is evolving beyond the Standard Model.
Stable configurations of classical (non-quantized)
nonlinear spinor fields
are also considered in quantum gravity \cite{krasnov2018gravity}.
(We note that nonzero Coulomb charge of a spinor field
ruins bi-frequency modes:
the charge--current density
$\bar\varPsi_{\pm\omega}\gamma^\mu\varPsi_{\pm\omega}$
of a bi-frequency mode \eqref{bf-0}
is time-dependent
-- unlike the scalar quantity $\bar\varPsi_{\pm\omega}\varPsi_{\pm\omega}$
-- and would radiate the energy via
electromagnetic field.)

%At the same time,
%the stability properties 
%require careful attention,
%% While the linear (``spectral'') stability
%% of one-frequency modes
%% can be established in certain cases
%% \cite{PhysRevLett.116.214101},
%%berkolaiko2015vakhitov,linear-b},

In the present article,
we are going to
(1) develop an approach to the (linear) stability
of one-frequency solitary waves in the NLD in (3+1)D;
(2)
present the numerical results
which show the linear stability of the NLD one-frequency nonlinear modes
and consequently the linear stability for
the DKG modes
for a wide range of parameters;
(3)
show that these stability results
imply (linear) stability of bi-frequency modes.

Let us emphasize that conclusive results
on the linear stability of localized modes in the NLD
have only been available in spatial dimensions $1$ and $2$
\cite{berkolaiko2015vakhitov,PhysRevLett.116.214101}.
In dimension three (and higher),
results on linear stability of solitary waves
(neither for one-frequency nor for bi-frequency ones)
were not available,
except for small amplitude one-frequency solitary waves
(``the nonrelativistic limit $\omega\lesssim m$'')
\cite{comech2014linear,boussaid2019spectral}
(and, as the matter of fact,
in this limit the cubic NLD is linearly unstable).

The approach to linear stability
of bi-frequency modes has been absent,
their stability properties were not known.
%(in spite of an earlier attempt
%\cite{boussaid2018spectral}).

\section{Nonlinear Dirac equation and Dirac--Klein--Gordon
system}
\label{sect-dkg}

The nonlinear Dirac equation
with cubic nonlinearity
(``the Soler model'' \cite{jetp.8.260,PhysRevD.1.2766})
is described by the Lagrangian density
\begin{align}\label{nld-density}
\mathscr{L}_{\mathrm{NLD}}
=
\bar\psi(\jj\gamma^\mu\p_\mu-M_s)\psi
+\frac{1}{2}(\bar\psi\psi)^2,
\qquad
\psi(t,x)\in\C^4,
\quad
(t,x)\in\R\times\R^3.
\end{align}
\begin{comment}
with the corresponding energy
\begin{align}
E_{\mathrm{NLD}}(\psi,\phi)=
\int_{\R^3}
\Big\{\psi^\dag(-\jj\bm\alpha\cdot\bm\nabla+M_s\beta)\psi
-\frac 1 2(\bar\psi\psi)^2
\Big\}\,dx.
\end{align}
\end{comment}
Here $\bar\psi=\psi^\dagger\beta$
is the Dirac conjugate of $\psi\in\C^4$,
with $\psi^\dagger$ denoting Hermitian conjugate of $\psi$.
We follow the standard convention
$\gamma^0=\beta=
\Big[
\begin{matrix}I_2&0\\[-4pt]0&-I_2\end{matrix}
\Big]$,
$\alpha^j=
\Big[
\begin{matrix}0&\sigma_j\\[-4pt]\sigma_j&0\end{matrix}
\Big]
$,
$\gamma^j=\beta\alpha^j$,
$1\le j\le 3$,
with $\sigma_j$ the Pauli matrices.
The NLD has important similarities
with the DKG system
described by the Yukawa-type
Lagrangian density
\begin{align}\label{dkg-density}
\mathscr{L}_{\mathrm{DKG}}
=
\bar\psi(\jj\gamma^\mu\p_\mu-M_s)\psi
+
\mbox{$\frac 1 2$}
\p^\mu\phi\,\p_\mu\phi
-
\mbox{$\frac 1 2$}
M_B^2\phi^2+g\phi\bar\psi\psi,
\end{align}
$\psi(t,x)\in\C^4$,
$\phi(t,x)\in\R$;
without loss of generality
we may assume that $g>0$
(the sign of $g$ does not affect the dynamics).
\begin{comment}
The corresponding energy is
\begin{align}
E_{\mathrm{DKG}}(\psi,\phi)=
\int_{\R^3}
\Big\{\psi^\dag(-\jj\bm\alpha\cdot\bm\nabla+M_s\beta)\psi
+\frac{\abs{\p_t\phi}^2
+\abs{\bm\nabla\phi}^2
+M_B^2\abs{\phi}^2}{2}
+g\phi\bar\psi\psi
\Big\}\,dx.
\end{align}
\end{comment}
Here are some of the properties
shared by both the NLD and the DKG models
\eqref{nld-density}
and
\eqref{dkg-density}:

\begin{enumerate}
\itemsep-5pt
%\smallskip
%\noindent{\it 1.\ }
\item
Both models have $\mathbf{SU}(1,1)$ symmetry
\eqref{galindo}.
By Noether's theorem, this continuous symmetry
leads to the conservation
of the $\mathbf{U}(1)$-charge
and the complex-valued ``Soler charge'',
%and of a complex-valued quantity
\begin{align}\label{charges}
Q=\int_{\R^3}\psi^\dag\psi\,dx,
\qquad
\varSigma=\int_{\R^3}
%\psi^\dag\jj\gamma^2\bmK\psi\,dx,
(\jj\gamma^2\bmK\psi)^\dag\psi\,dx,
\end{align}
with $\bmK$ the complex conjugation.
Under the
$\mathbf{SU}(1,1)$ symmetry transformation
\eqref{galindo},
$\psi\mapsto\mathbf{g}\psi=(a+\jj b\gamma^2\bmK)\psi$,
the energy is conserved,
while
the charges \eqref{charges}
of $\psi$ and $\mathbf{g}\psi$
are related by
\begin{align*}
Q(\mathbf{g}\psi)
=(\abs{a}^2+\abs{b}^2)Q(\psi)
+2\Re(a\bar b \varSigma(\psi)),
\qquad
\varSigma(\mathbf{g}\psi)
=a^2\varSigma(\psi)
+2a b Q(\psi)
+b^2\overline{\varSigma(\psi)}.
\end{align*}
In particular, one has
$Q(\mathbf{g}\psi)^2-\abs{\varSigma(\mathbf{g}\psi)}^2
=
Q(\psi)^2-\abs{\varSigma(\psi)}^2$
(see \cite{boussaid2018spectral}).

%\smallskip
%\noindent{\it 2.\ }
\item
%Besides one-frequency modes
Both the DKG system and the NLD
have bi-frequency modes
(see \eqref{bf-1} below).
%with one-frequency modes
%(of the form \eqref{sw} below)
%being their particular case.
Generically,
bi-frequency solitary waves
are not of the form $e^{\Omega t}\upvarphi(x)$,
with $\Omega$ a Lie algebra element;
while both models
are $\mathbf{SU}(1,1)$-invariant,
the set of
\mbox{localized solutions admits}
a larger symmetry group $\mathbf{SU}(2,2)$
(this corresponds
to the choice of
$\bm\xi,\,\bm\eta\in\C^2$ in \eqref{bf-1}
satisfying $\norm{\bm\xi}^2-\norm{\bm\eta}^2=1$).
One-frequency modes
$\psi_{\omega}(t,x)$
appear as the ``endpoints''
of the manifold of bi-frequency modes
$\varPsi_{\pm\omega}(t,x)$
(\eqref{bf-1} with $\bm\eta=0$ turns into \eqref{sw}).

%\smallskip
%\noindent{\it 3.\ }
\item
Given the linear stability
of one- and bi-frequency solitary waves in
the NLD for
% a wide range of parameters,
$\omega\in(0.254M_s,0.936M_s)$
(see numerics in Section~\ref{sect-numerics}),
by the perturbation theory,
there is a nearby stability region
for one- and bi-frequency solitary waves
in the DKG system
for these values of $\omega$
and for $\abs{g}\approx M_B$ sufficiently large
(and possibly for some other values
of $g$, $M_B$, and $\omega$).
\end{enumerate}

\vskip -2.5mm

%Let us now make several comments
%about the DKG system
%with the Lagrangian density \eqref{dkg-density},
%where without loss of generality
%we may assume that $g>0$
%(the sign of $g$ does not affect the dynamics).
The mass spectrum for the NLD is given in
\cite{PhysRevD.1.2766}
(the energy has a minimum at
$\omega\approx 0.936 M_s$ and grows arbitrarily large
as $\omega\to 0$ or $\omega\to M_s$).
To the best of our knowledge,
the energy spectrum for the DKG has not been computed
yet; for
$\omega\to M_s$,
the asymptotic behavior 
would parallel that of the NLD
%give the answer similar to Dirac--Maxwell system
%(similar to \cite{comech2018small},
%(similar to \cite{boussaid2019spectral})
(see \cite{boussaid2017nonrelativistic}),
yielding
$E(\omega),\,Q(\omega)\propto(M_s-\omega)^{-1/2}$.
%for $\omega<M_s$, $\omega\to M_s$.
For the weakly relativistic
%(one-frequency)
solitary waves
$\psi(t,x)=e^{-\jj\omega t}\upvarphi(x)$
such that
$\omega\lesssim M_s$,
$\abs{\upvarphi(x)}$ is small,
and $\bar\upvarphi(x)\upvarphi(x)>0$
for all $x\in\R^3$,
the equation
$(-\Delta+M_B^2)\phi=g\bar\psi\psi$
shows that
%let us mention
%that weakly relativistic solitary waves
%the latter quantity is sign-definite)
the (time-independent) boson field $\phi(x)$ is positive.
%%AC ????? 
Since it is this field that leads to the formation
of a localized state, we conclude that
(at least in the limit of small amplitudes)
the particles which generate the boson field $\phi$
of the same sign
would be attracting
(and one can show that there are no localized solitary waves
$\psi(t,x)=\upvarphi(x)e^{-\mathrm{i}\omega t}$
with $\bar\upvarphi(x)\upvarphi(x)<0$ for all $x$).
The above reasoning also applies to bi-frequency solutions:
bi-frequency solitary waves $\psi_1$ and $\psi_2$
of the form \eqref{bf-1}
are mutually attracting when $\bar\psi_1\psi_1$
\mbox{and $\bar\psi_2\psi_2$ are strictly positive.}

Since the construction of bi-frequency modes
in the DKG system and in the NLD
is the same,
% while the latter system is shorter to describe,
in the remainder of the article
we discuss the linear stability
concentrating on the NLD.
We point out that the DKG system
turns into the NLD
in the limit of heavy
% Klein--Gordon
bosons
and large coupling constants,
$|g|\approx M_B\to\infty$,
when the interaction term
$g\phi\psi
\sim g((\p_t^2-\Delta+M_B^2)^{-1}g\bar\psi\psi)\beta\psi$
in the equation for $\psi$
%in the Lagrangian
turns into the scalar-type self-interaction
term $\sim (g^2\bar\psi\psi/M_B^2)\beta\psi$
in the NLD.
%% \[
%% \bar\psi(\jj\gamma^\mu\p_\mu-M_s)\psi
%% +
%% \mbox{$\frac 1 2$}
%% \p^\mu\phi\,\p_\mu\phi
%% -
%% \mbox{$\frac 1 2$}
%% M_B^2\phi^2
%% -g\phi\bar\psi\psi;
%% \]
%$\psi^\dag D_0\psi$
In this limit,
the shape of localized
spinor modes of the DKG approaches
that in the NLD;
the same convergence takes place
for the operators corresponding to the linearization
at a localized mode and
hence for the linear stability properties.
The approximation of the DKG system
with the NLD is justified
if the mass $M_s$ of the spinor field
is much smaller than the mass $M_B$
of the Klein--Gordon field,
with the coupling constant $g\sim M_B$.
For example,
this would be justified for $M_s$
just above the
Lee--Weinberg lower bound of
$\sim 2\,\mathrm{GeV}$
%for the mass of stable neutrinos
%(or perhaps up to $\sim 13\,\mathrm{GeV}$)
for the DM neutrinos,
or perhaps from $1.3$ to $13\,\mathrm{GeV}$
\cite{kolb1986lee,aalbers2022next},
%\cite{kolb1986lee}
%\cite{kolb1986lee,aalbers2022next}
while $M_B$ corresponds to the Higgs boson at $125\,\mathrm{GeV}$.
%\end{comment}

%\newpage

\section{Linear stability of one-frequency spinor modes}

We consider the cubic NLD
\cite{jetp.8.260,PhysRevD.1.2766}
\begin{equation}\label{Soler}
\jj\p\sb{t}\psi
=-\jj\bm\alpha\cdot\bm\nabla\psi+
%f(\bar\psi\psi)
M_s\beta\psi
-(\bar\psi\psi)\beta\psi,
\qquad
\psi(t,x)\in\C^4,
\end{equation}
with
$M_s>0$ the mass of the spinor field.
There are solitary wave solutions
to \eqref{Soler}
%localized modes
%$\upvarphi(x)e^{-\jj\omega t}$,
%with
of the form
\cite{PhysRevD.1.2766}
\begin{align}\label{sw}
\psi_{\omega}(t,x)
=
\upvarphi(x)e^{-\jj\omega t},
\qquad
\upvarphi(x)
=
\begin{bmatrix}
v(r,\omega)\bm\xi
\\
\jj\sigma_r u(r,\omega)\bm\xi
\end{bmatrix}
e^{-\jj\omega t}
,
\end{align}
where
$r=\abs{x}$,
$\sigma_r=r^{-1}\bm{x}\cdot\bm\sigma$,
$\bm\xi\in\C^2$,
$\norm{\bm\xi}=1$;
the scalar functions
$v(r,\omega),\,u(r,\omega)$,
are real-valued
%even and odd, respectively,
and satisfy
(cf. \cite{esteban1995stationary,boussaid2017nonrelativistic})
%% \[
%% \begin{cases}
%% \omega v=\p_r u+\frac{2}{r}u+f(v^2-u^2)v,
%% \\
%% \omega u=-\p_r v-f(v^2-u^2)u,
%% \end{cases}
%% \qquad
%% r>0
%% \]
\begin{align}\label{vu-such}
%\omega v=\p_r u+\frac{2u}{r}+f(v^2-u^2)v,
%\omega v=\p_r u+2r^{-1}u+f(\bar\upvarphi\upvarphi)v,
\omega v=\p_r u+2r^{-1}u+(M_s-(v^2-u^2))v,
\qquad
%\omega u=-\p_r v-f(\bar\upvarphi\upvarphi)u.
\omega u=-\p_r v-(M_s-(v^2-u^2))u.
\end{align}
%where $\bar\upvarphi\upvarphi=v^2-u^2$.
%\subsection{Bi-frequency modes}
%In this section, we assume that the spatial dimension is $n=3$.
\begin{comment}
A bi-frequency mode
can be obtained
from a one-frequency mode
$
e^{-\jj\omega t}
\begin{bmatrix}
v(r)\bm\xi_0
\\
\jj u(r)\sigma_r\bm\xi_0
\end{bmatrix}$,
$\bm\xi_0\in\C^2$,
$\norm{\bm\xi_0}=1$,
by the application of the $\mathbf{SU}(1,1)$ transformation:
given $a,\,b\in\C$,
$\abs{a}^2-\abs{b}^2=1$,
one has
\begin{equation}\label{bf}
\Big(
a-b
\begin{bmatrix}0&\!\!\sigma_2\\-\sigma_2&\!\!0\end{bmatrix}
\bmK\Big)
e^{-\jj\omega t}
\begin{bmatrix}
v(r)\bm\xi_0
\\
\jj u(r)\sigma_r\bm\xi_0
\end{bmatrix}
=
e^{-\jj\omega t}
\begin{bmatrix}v(r)a\bm\xi_0
\\
\jj u(r)\sigma_r a\bm\xi_0
\end{bmatrix}
+
e^{\jj\omega t}
\begin{bmatrix}
-\jj u(r)\sigma_r b\sigma_2\bmK\bm\xi_0
\\
v(r)b\sigma_2\bmK\bm\xi_0
\end{bmatrix}.
\end{equation}
(We recall that $\sigma_2\bmK$ commutes with
$\jj\sigma_j$, $1\le j\le 3$.)
\end{comment}
\begin{comment}
Similarly,
there are solitary waves of the form
\begin{align}\label{sw-2}
\upchi(x)e^{\jj\omega t},
\qquad
\upchi(x)
=\begin{bmatrix}
-\jj\sigma_r u(r,\omega)\bm\eta
\\
v(r)\bm\eta
\end{bmatrix},
\end{align}
with $\bm\eta\in\C^2$, $\norm{\bm\eta}=1$,
and with $v,\,u$
satisfying equations similar to \eqref{vu-such};
the stability analysis of 
\eqref{sw}
and \eqref{sw-2}
is the same so we always consider
the (one-frequency) solitary waves
of the first type.
\end{comment}

We recall the
% usual
%Bogoliubov
%--de Gennes
linear stability analysis
of standard, one-frequency modes:
%\cite{zakharov-1967,kolokolov-1973}.
given a solitary wave
$e^{-\jj\omega t}\upvarphi$
(or,
more generally,
$e^{\Omega t}(\upvarphi+\rho(t))$,
with $\Omega$
from the Lie algebra
of the symmetry group $G$ of the Lagrangian),
one considers its perturbation in the form
$(\upvarphi+\rho(t))e^{-\jj\omega t}$
(or,
more generally,
$e^{\Omega t}(\upvarphi+\rho(t))$),
%,
%with $\Omega$
%from the Lie algebra
%of the symmetry group of the Lagrangian),
writes a linearized equation on $\rho$,
\ak{$
  \p_t
  \begin{bmatrix}\rho\\\bar\rho\end{bmatrix}
  =
  A\begin{bmatrix}\rho\\\bar\rho\end{bmatrix},
  $
  with $\bar\rho\,$ the complex conjugate of $\rho$
  (so that the operator $A$ is $\C$-linear),
}
and studies the spectrum of
the corresponding operator
\ak{$A$}
(which does not depend on $t$
due to the
% $\mathbf{U}(1)$-invariance
$G$-invariance
of the original system).
If the spectrum
\ak{of $A$}
is purely imaginary,
one says that the solitary wave is \emph{spectrally stable}
(or \emph{linearly stable}).
%which is the weakest form of stability.
%%\end{comment}
\begin{comment}
For a bi-frequency mode
\eqref{bf-0},
a natural idea is to write the solution
in the form
$
(\upvarphi+\rho_1(t))
e^{-\jj\omega t}+(\upchi+\rho_2(t)e^{\jj\omega t}$;
%% writing the linearized equation
%% \[
%% \p_t
%% \begin{bmatrix}R\\\bar R\end{bmatrix}
%% =
%% A\begin{bmatrix}R\\\bar R\end{bmatrix},
%% \qquad
%% \mbox{with}
%% \quad
%% R=\begin{bmatrix}\rho_1\\\rho_2\end{bmatrix},
%% \]
%% and again study the spectrum of $A$.
the difficulty is
in choosing the decomposition of the
total perturbation $\rho$ into
$\rho(t)
=\rho_1(t)e^{-\jj\omega t}+\rho_2(t)e^{\jj\omega t}$
%initial perturbation $\rho_0(x)$
%into the sum $\rho_0(x)=\rho_1(0,x)+\rho_2(0,x)$
in such a way that
the linearized equation
%$A$
does not contain time dependence
in the form of factors $e^{\pm\jj\omega t}$.
\end{comment}
\ak{\sout{
In this article, we achieve this
in the framework of interacting spinor and scalar fields
(and similarly in the case of spinor fields with
self-interaction of scalar type)
and then show that
in certain cases the stability indeed takes place.
}}
%\subsection{Stability of one-frequency modes via radial reduction}
Consider a perturbation
of a one-frequency solitary wave
\eqref{sw},
$\big(
\upvarphi(x)+
\rho(t,x)
\big)e^{-\jj\omega t}$,
$\rho(t,x)\in\C^4$.
The linearization at $\upvarphi e^{-\jj\omega t}$
-- that is, 
the linearized equation on
$\rho$
--
takes the form
\begin{align*}
%\label{lin}
\jj\p_t\rho
\!=\!
\calL\rho
\!:=D_0\rho
+%f(\bar\upvarphi\upvarphi)
(M_s-\bar\upvarphi\upvarphi)
\beta\rho
%+2f'(\bar\upvarphi\upvarphi)
-2
\beta\upvarphi\Re(\bar\upvarphi\rho)
-\omega\rho.
%\qquad
%\dom(\calL)=H^1(\R^3,\C^4),
\end{align*}
%with
%$f=f(\bar\upvarphi\upvarphi)$,
%$f'=f'(\bar\upvarphi\upvarphi)$.
Note that the operator
$\calL$
%in \eqref{lin}
is not $\C$-linear
because of the term $\Re(\bar\upvarphi\rho)$.
\begin{comment}
Besides
$\calL$,
we also introduce a
$\C$-linear
selfadjoint operator
\begin{align}\label{def-l0}
\calL_0=D_0+f\beta-\omega.
%\qquad\dom(\calL_0)=H^1(\R^3,\C^4).
\end{align}
\end{comment}
It turns out that
$\calL$ has the following
invariant subspaces
for $\ell\ge 0$, $-\ell\le m\le\ell$:
%\begin{widetext}
\begin{align}
\label{def-x}
\scrX_{\ell,m}
&=
\left\{
\sum\limits\sb{\pm}
\begin{bmatrix}
(a_{\pm m}+p_{\pm m}\GG)\hh_{\ell}^{\pm m}\e_1
\\
\jj\sigma_r(b_{\pm m}+q_{\pm m}\GG)\hh_{\ell}^{\pm m}\e_1
\end{bmatrix}
\right\},
\qquad
%\\
%\label{def-y}
\scrY_{\ell}
%&
=
\Biggl\{
\begin{bmatrix}
R\hh_{\ell}^{-\ell}\e_2
\\
\jj \sigma_r S\hh_{\ell}^{-\ell}\e_2
\end{bmatrix}
\Biggr\}.
\end{align}
%\end{widetext}
%are invariant subspaces of $\calL$.
Above,
$\GG$
is the angular part
of $\bm\sigma\cdot\bm\nabla$,
defined by the relation
\begin{align}
\bm\sigma\cdot\bm\nabla
=\sigma_r\Big(\p_r-\frac{\GG}{r}\Big),
\qquad
\sigma_r=r^{-1}\bm{x}\cdot\bm\sigma;
\end{align}
$Y_{\ell}^m
=\sqrt{\frac{(2\ell+1)(\ell-|m|)!}{4\pi(\ell+|m|)!}}
e^{\jj m\phi}P_\ell^{|m|}(\cos\theta)$
are spherical harmonics
of degree $\ell\ge 0$
and order $\abs{m}\le\ell$
(with $P_\ell^m$ the associated Legendre polynomials);
$a_{\pm m},\,\dots,\,R,\,S$
are functions of $r$.
We note that
$\GG$
is related to the operator of spin-orbit interaction by
$
2\bm{S}\cdot\bm{L}
=
\begin{bmatrix}
\GG&0\\0&\GG
\end{bmatrix}
$,
with
$\bm{S}=-\frac{\jj}{4}\bm\alpha\wedge\bm\alpha$
the spin angular momentum operator
and
$\bm{L}=\bm{x}\wedge(-\jj\bm\nabla)$
the orbital angular momentum operator
%$G=-\sum_{1\le j<k\le 3}\jj\alpha^j\alpha^k(x_j p_k-x_k p_j)$
%the spin-orbit operator
\cite{thaller1992dirac}
%(for the properties of $\GG$, see
%also
% \cite{decarli1999strong}
(see also \cite{kalf1999note}).
While all the invariant spaces
$\scrX_{\ell,m}$,
$\scrY_{\ell}$
are needed
to represent an arbitrary perturbation
of a solitary wave,
$\scrY_\ell$
can be discarded
\newpage
 from future consideration:
the restriction of $\calL$
onto $\scrY_\ell$ coincides with
selfadjoint operator
\begin{align}
\label{def-l0}
\calL_0=D_0
%+f(\bar\upvarphi\upvarphi)
+(M_s-\bar\upvarphi\upvarphi)
\beta-\omega,
\end{align}
hence the equation
$\jj\p_t\rho=\calL\rho$ restricted onto
$\scrY_\ell$
does not have modes growing exponentially in time
so cannot lead to linear instability.
%$\sigma\big(-\jj\calL\at{\scrY_\ell}\big)\subset\jj\R$.

%The operator $\calL_0(\omega)$
%is selfadjoint and $\C$-linear.
In the space
$\scrX_{\ell,0}$,
$\ell\ge 1$,
acting on vectors
$\Psi=(a_0,b_0,p_0,q_0)^T$
with components
%(the superscript $T$ denotes the transpose)
depending on $t$ and $r$,
the operator $\calL_0(\omega)$
is represented by the matrix-valued operator
%\nb{$L_0(\omega,\ell)$ or $\calL_0(\omega)$?}
\begin{align}\label{def-l}
L_0(\omega,\ell)
=
\begin{bmatrix}
f-\omega&\p_r+\frac{2}{r}&0&\frac{\ell(\ell+1)}{r}
\\
-\p_r&-f-\omega&\frac{\ell(\ell+1)}{r}&0
\\
0&\frac{1}{r}&f-\omega&\p_r+\frac{1}{r}
\\
\frac{1}{r}&0&-\p_r-\frac{1}{r}&-f-\omega
\end{bmatrix},
\end{align}
where
$f
%=f(\bar\upvarphi\upvarphi)
=M_s-\bar\upvarphi\upvarphi
$.
Since
$\calL(\omega)$ is not $\C$-linear,
we introduce the
$\C$-linear operator
$\frakL(\omega)$ such that
$
\begin{bmatrix}
\calL\rho
\\
\bmK\calL\rho
\end{bmatrix}
=\frakL(\omega)
\begin{bmatrix}
\rho\\
\bmK\rho
\end{bmatrix}
$.
%Similarly, one can consider the spaces
%$\scrX_{\ell,m}$ with $m\ne 0$.
%the operator 
%(on vectors with components $a_{\pm m}(t,r),\,\dots,\,q_{\pm m}(t,r)$,
%respectively),
Perturbations corresponding to spherical
harmonics of degree $\ell$
and orders $\pm m$ are mixed:
the linearized equation contains
$\Psi_{\ell,m}$
and
$\bmK\Psi_{\ell,-m}$.
\ak{\sout{It turns out that the
 term $\Re(\bar\upvarphi\rho)$
 in \eqref{lin}
 results in the mixing of
 $\Psi_{\pm m}$ and $\bmK\Psi_{\mp m}$
 but not of e.g.
 $\Psi_{m}$ and $\bmK\Psi_{m}$.
}}
%% allowing to reduce the size
%% of matrices under consideration
%% from $16\times 16$
%% to $8\times 8$.
%the equation $\jj\rho=\calL\rho$
%In the space
%$1\le m\le \ell$,
When acting on vectors
$\begin{bmatrix}\Psi_{m}\\\bmK\Psi_{-m}\end{bmatrix}
%\in\scrX_{\ell,m}
$,
with
$
\Psi_{m}=(a_{m},b_{m},p_{m},q_{m})^T$,
$\frakL(\omega)$
is represented by
\begin{align}\label{def-w}
&
\begin{bmatrix}
L_0(\omega,\ell)&0\\[1ex]
0&L_0(\omega,\ell)
\end{bmatrix}
+
\begin{bmatrix}
V&m V&V&-m V
\\
0&0&0&0
\\
m V&V&-m V&V
\\
0&0&0&0
\end{bmatrix},
%\nonumber
%\\
%&%with $V$ given by
\qquad
V(r,\omega):=
-\begin{bmatrix}v^2&-uv\\-uv&u^2\end{bmatrix},
%% \quad
%% v=v(r,\omega),
%% \quad
%% u=u(r,\omega),
\end{align}
with $L_0$ from \eqref{def-l}
and
with $v,\,u$
corresponding to the profile of the solitary wave
\eqref{sw}.
The linear stability of
$\jj\p_t\rho=\calL\rho$
reduces to studying the linear stability
in each of the invariant subspaces
$\scrX_{\ell,m}$,
$\ell\ge 1$, $-\ell\le m\le\ell$,
which in turn reduces to studying the spectrum
of operators $\mathbf{A}_{\ell,m}$
given by
\begin{align}\label{lo}
-\jj\left\{
\begin{bmatrix}L_0&0\\[1ex]0&-L_0\end{bmatrix}
+
\begin{bmatrix}
V&m V&V&-m V
\\
0&0&0&0
\\
-V&-m V&-V&m V
\\
0&0&0&0
\end{bmatrix}
\right\},
\end{align}
with $L_0$ from \eqref{def-l}
and $V$ from \eqref{def-w}.
We note that the eigenvalues of $\mathbf{A}_{\ell,\pm m}$
are mutually complex conjugate.

The case $\ell=0$ is exceptional:
the corresponding perturbations
have the same angular structure
as the solitary wave itself
and allow a simpler treatment
\cite{PhysRevLett.116.214101}.
In that case, $Y_0^0=1$,
so in \eqref{def-x}
one takes $p_0(r)=q_0(r)=0$;
instead of $L_{0}$ from \eqref{def-l}
one needs to consider
\[
L_{00}(\omega)=
\begin{bmatrix}
f-\omega&\p_r+\frac{2}{r}\\
-\p_r&-f-\omega
\end{bmatrix},
\qquad
f=M_s-\bar\upvarphi\upvarphi,
\]
and for the linear stability
with respect to perturbations from
$\scrX_{0,0}$
one needs to study the spectrum of
(cf. \eqref{lo})
\[
\mathbf{A}_{00}
=
-\jj
\begin{bmatrix}
L_{00}+V&V
\\
-V&-L_{00}-V
\end{bmatrix}.
\]
%which is of smaller size than \eqref{lo}.

\section{Linear stability of bi-frequency spinor modes}

By \cite{boussaid2018spectral},
if \eqref{sw} is a solitary wave solution to
\eqref{Soler}, then so is
a bi-frequency solitary wave
or \emph{bi-frequency mode},
\begin{align}\label{bf-1}
\varPsi_{\pm\omega}(t,x)
=
\begin{bmatrix}
v(r)\bm\xi
\\
\jj u(r)\sigma_r\bm\xi
\end{bmatrix}
e^{-\jj\omega t}
+
\begin{bmatrix}
-\jj u(r)\sigma_r\bm\eta
\\
v(r)\bm\eta
\end{bmatrix}
e^{\jj\omega t},
\end{align}
with
$\bm\xi,\,\bm\eta\in\C^2$,
$\norm{\bm\xi}^2-\norm{\bm\eta}^2=1$.
If $\bm\xi,\,\bm\eta\in\C^2$
in the expression for the bi-frequency mode \eqref{bf-1}
are mutually orthogonal,
then this mode can be obtained from one-frequency mode \eqref{sw}
by application of the transformation
from the symmetry group $\mathbf{SU}(1,1)$ of the NLD
(see \cite{galindo1977remarkable,boussaid2019nonlinear}),
given by
\eqref{galindo}.
%% \begin{align}\label{su11}
%% a+b
%% \begin{bmatrix}0&\!\!\sigma_2\\-\sigma_2&\!\!0\end{bmatrix}
%% \bmK,
%% \quad
%% a,\,b\in\C,
%% \quad
%% \abs{a}^2-\abs{b}^2=1.
%% \end{align}
In this case, the stability of \eqref{bf-1}
follows from the corresponding result for
\eqref{sw}
by applying to a perturbed bi-frequency mode
the $\mathbf{SU}(1,1)$ symmetry transformation \eqref{galindo}
which makes one-frequency solution from
a bi-frequency one.
If $\langle\bm\xi,\bm\eta\rangle\ne 0$,
though,
then a bi-frequency solitary wave \eqref{bf-1}
can not be obtained from \eqref{sw}
via the action of
$\mathbf{SU}(1,1)$;
in this case, stability analysis of \eqref{bf-1}
does not reduce to the stability analysis of \eqref{sw}.
It turns out, though, that
the symmetry transformation \eqref{galindo}
can be used to reduce a bi-frequency solitary wave
\eqref{bf-1}
to the case when $\bm\xi$ and $\bm\eta$ are parallel;
thus, to study the linear stability of bi-frequency
solitary waves,
it is enough to concentrate on the two cases:
when $\bm\xi$ and $\bm\eta$ are mutually orthogonal
(the case equivalent to one-frequency solitary waves)
and when
$\bm\xi$ and $\bm\eta$ are parallel
(a nontrivial case).

As we pointed out,
bi-frequency solitary waves \eqref{bf-1}
are not of the form $e^{\Omega t}\upvarphi(x)$,
with $\upvarphi(x)$ a stationary solution
and $\Omega$ the element of the Lie algebra of the
symmetry group of the system,
so one cannot use the Grillakis--Shatah--Strauss approach
by considering the Ansatz
$e^{\Omega t}\big(\upvarphi(x)+\varrho(t,x)\big)$
and studying the linearized equation on the perturbation $\varrho$;
here, instead, we develop an ad hoc approach.
We consider a perturbation of a bi-frequency mode
\eqref{bf-1}
in the form $\psi(t,x)=\varPsi_{\pm\omega}(t,x)+\varrho(t,x)$,
%% =
%% \Big(
%% \begin{bmatrix}
%% v(r)\bm\xi
%% \\
%% \jj u(r)\sigma_r\bm\xi
%% \end{bmatrix}
%% +\rho_1(t,x)
%% \Big)
%% e^{-\jj\omega t}
%% +
%% \Big(
%% \begin{bmatrix}
%% -\jj u(r)\sigma_r\bm\eta
%% \\
%% v(r)\bm\eta
%% \end{bmatrix}
%% +\rho_2(t,x)
%% \Big)
%% e^{\jj\omega t}.
%The perturbation satisfies
% the linearized equation
where $\varrho(t,x)$ satisfies
%(cf. \eqref{lin})
\begin{align}\label{lin-2}
\jj\p_t\varrho
=D_0\varrho
%+f(\bar\varPsi_{\pm\omega}\varPsi_{\pm\omega})
+(M_s-\bar\varPsi_{\pm\omega}\varPsi_{\pm\omega})
\beta\varrho
%+2f'(\bar\varPsi_{\pm\omega}\varPsi_{\pm\omega})
-2
\Re(\bar\varPsi_{\pm\omega}\varrho)
\beta\varPsi_{\pm\omega}.
\end{align}
%where $f=f(\bar\varPsi_{\pm\omega}\varPsi_{\pm\omega})$ and $f'=f'(\bar\varPsi_{\pm\omega}\varPsi_{\pm\omega})$.
We note that
for $\varPsi_{\pm\omega}$ from \eqref{bf-1},
$\bar\varPsi_{\pm\omega}\varPsi_{\pm\omega}$ does not depend on time.
For each $\ell\in\N_0$,
the linearization \eqref{lin-2}
is invariant in the spaces
formed by
$\varrho_1(t,x)$
and by
$\varrho_1(t,x)+\varrho_2(t,x)+\varrho_3(t,x)$,
where
%\begin{widetext}
\begin{align}
%\label{def-rho-1-2}
\label{def-rho-1}
\varrho_1(t,x)
&=
\sum_{m=-\ell}^\ell
\left\{
\begin{bmatrix}
(a_{m}+p_{m}\GG)\hh_{\ell}^{m}\bm\xi
\\
\jj\sigma_r(b_{m}+q_{m}\GG)\hh_{\ell}^{m}\bm\xi
\end{bmatrix}
e^{-\jj\omega t}
+
\begin{bmatrix}
-\jj\sigma_r(\bar b_{m}+\bar q_{m}\GG)
\hh_{\ell}^{-m}\bm\eta
\\
(\bar a_{m}+\bar p_{m}\GG)\hh_{\ell}^{-m}\bm\eta
\end{bmatrix}
e^{\jj\omega t}
\right\},
%\quad\mbox{and}\quad
\\
\label{def-rho-2}
\varrho_2(t,x)
&=
\begin{bmatrix}
R
\hh_{\ell}^{-\ell}\bm\xi^\perp
\\
\jj\sigma_r
S
\hh_{\ell}^{-\ell}\bm\xi^\perp
\end{bmatrix}
e^{\jj\omega t}
,
\qquad
\varrho_3(t,x)
=
\begin{bmatrix}
-\jj\sigma_r U
%% \hh_{\ell}^{\ell}\bm\eta^\perp
%% \\
%% T
%% \hh_{\ell}^{\ell}\bm\eta^\perp
 \tilde\hh_{\ell}^{-\ell}\bm\eta^\perp
 \\
 T
 \tilde\hh_{\ell}^{-\ell}\bm\eta^\perp
\end{bmatrix}
e^{-\jj\omega t}
,
\end{align}
%\end{widetext}
with
$a_{m}$, $b_{m}$, $p_{m}$, $q_{m}$
(with $\abs{m}\le\ell$),
%$\ell\in\N_0$,
$R$, $S$, $T$, and $U$
complex-valued 
functions of $t$ and $r$
and with $\bm\xi,\,\bm\eta\in\C^2$ from \eqref{bf-1}.
It is assumed that $\bm\xi$ is parallel to $\bm{e}_1$,
while the notation $\tilde\hh_{\ell}^{-\ell}$
refers to the $(\ell,-\ell)$-spherical harmonic
in the rotated coordinate system in $\mathbb{R}^3$
such that the corresponding transformation of $\mathbb{C}^2$
makes $\bm\eta$ parallel to $\bm{e}_1$.
We point out that
any perturbation can be decomposed into
$\varrho_1+\varrho_2+\varrho_3$
summed over $\ell\ge 0$.
The invariance in these subspaces
is to be understood in the sense that
there is a time-independent,
\nopagebreak
$\R$-linear (but not $\C$-linear)
differential operator $\calA(x,\bm\nabla)$
such that equation \eqref{lin-2}
for $\varrho$
is equivalent to
$\p_t\Psi=\calA\Psi$,
where
$\Psi$ contains all of $a_m$, $b_m$, $\dots$,
with $\abs{m}\le\ell$.
(The expressions
%\eqref{def-rho-1-2}
\eqref{def-rho-1},
\eqref{def-rho-2}
are such that
$\Re(\bar\varPsi_{\pm\omega}\varrho)$
does not contain factors of $e^{\pm 2\jj\omega t}$,
so that \eqref{lin-2}
only contains two groups of terms,
with factors $e^{-\jj\omega t}$ and $e^{\jj\omega t}$.)
Moreover,
% we can assume that 
it is enough to consider perturbations
with
$\varrho_2=0$ and $\varrho_3=0$:
indeed, if $\varrho=\varrho_1+\varrho_2+\varrho_3$
%with e.g. $\varrho_2\ne 0$
satisfies
$\lambda\varrho=\calA\varrho$,
one can deduce that
%$\lambda\varrho_2=-\jj\calL_0\varrho_2$
$\lambda\varrho_2=-\jj(\calL_0+2\omega)\varrho_2$
and
$\lambda\varrho_3=-\jj\calL_0\varrho_3$,
with
$\calL_0$ from \eqref{def-l0}
being selfadjoint,
so the assumption
that either $\varrho_2\ne 0$
or $\varrho_3\ne 0$
leads to $\lambda\in\jj\R$,
causing no linear instability.
%similarly, we can assume that
%$\varrho_3=0$.
With 
$\varrho_2=0$
and $\varrho_3=0$,
we have:
\begin{align}\label{re-big}
\Re(\bar\varPsi_{\pm\omega}\varrho)
=
{\sum}_{\abs{m}\le\ell}
\Re\big[
(v a_{m} - u b_{m})\hh_{\ell}^{m}
+
(v p_{m}-u q_{m})\big(
\bm\xi^\dagger\GG\hh_{\ell}^{m}\bm\xi
+\bm\eta^\dagger\GG\hh_{\ell}^{m}\bm\eta
\big)
\big].
\end{align}
%% There is a differential operator $A(x,\bm\nabla)$
%% which is $\C$-linear such that
%% $\p_t\Psi=\calA\Psi$ can be written as
%% \[
%% \p_t
%% \begin{bmatrix}\Psi\\\bmK\Psi\end{bmatrix}
%% =A\begin{bmatrix}\Psi\\\bmK\Psi\end{bmatrix}.
%% \]
%% The spectrum of the (complexification)
%% of the operator $\calA$ of linearization at $\varPsi_{\pm\omega}$
%% In particular, there is no ``frequency mixing'',
%% in the sense that
%% the quantity $\bar\psi\psi$
%% in the first order
%% does not contain terms
%% with factors $e^{\pm 2\jj\omega t}$.
One can see from \eqref{re-big}
that the linear stability of one-frequency and bi-frequency modes
from perturbations
corresponding to spherical harmonics of degree zero
(same angular structure as the solitary wave itself)
is the same: $\GG Y_0^0=0$, hence the terms with
$\bm\xi$ and $\bm\eta$-dependence drop out.
%% We postpone the detailed analysis
%% of \eqref{re-big} for general $\bm\xi$ and $\bm\eta$
%% for future work,

Let us now consider harmonics of degree
$\ell\ge 1$.
As we already pointed out,
it is enough to focus on the two ``endpoint'' cases:
when $\bm\eta$ is parallel to $\jj\sigma_2\bmK\bm\xi$
(hence orthogonal to $\bm\xi$; this case can be
transformed via the $\mathbf{SU}(1,1)$ symmetry transformation 
to one-frequency solitary waves \eqref{sw})
and when $\bm\eta$ is parallel to $\bm\xi$.
%(the $\mathbf{SU}(1,1)$ symmetry transformation
%\eqref{galindo} can be used to reduce to
%one of these two cases).
In the first case,
one has
\[
\Re\big(
\bm\xi^\dagger\GG\hh_{\ell}^{m}\bm\xi
+
\bm\eta^\dagger\GG\hh_{\ell}^{m}\bm\eta
\big)
=
\Re\Big(
\bm\xi^\dagger\GG\hh_{\ell}^{m}\bm\xi
-
\norm{\bm\eta}^2
\frac{\bm\xi^\dagger\GG\hh_{\ell}^{m}\bm\xi}{\norm{\bm\xi}^2}
\Big)
=
\frac{1}{\norm{\bm\xi}^{2}}
\Re\big[
\bm\xi^\dagger
\GG\hh_{\ell}^{m}
\bm\xi
\big]
\]
and then the linearized operator
coincides with the linearization
at a one-frequency solitary wave
(corresponding to the spherical harmonic of degree $\ell$
and order $m$,
with the ``polarization''
given by $\bm\xi_0=\bm\xi/\norm{\bm\xi}\in\C^2$
in place of $\bm\xi$).
Indeed,
in this case the bi-frequency solitary wave
can be obtained from a one-frequency solitary wave
via application to \eqref{bf-1}
of an appropriate $\mathbf{SU}(1,1)$
transformation \eqref{galindo}, hence
the one-frequency and bi-frequency modes
share their stability properties.
(We note that
if $\bm\xi_0=(1,0)^T$,
then
$\Re(\bm\xi_0^\dagger\GG\hh_{\ell}^{m}\bm\xi_0)=m$;
this is what leads to factors of $m$ in \eqref{lo}.)

If instead
$\bm\xi$ and $\bm\eta$
are parallel
(in this case, the bi-frequency solitary wave
\emph{cannot}
be obtained from a one-frequency solitary wave
with the aid of
the $\mathbf{SU}(1,1)$ transformation), then
the part with $\bm\xi$ and $\bm\eta$ from
\eqref{re-big}
takes the form
\[
\bm\xi^\dagger\GG\hh_{\ell}^{m}\bm\xi
+
\bm\eta^\dagger\GG\hh_{\ell}^{m}\bm\eta
=
\big(\norm{\bm\xi}^2+\norm{\bm\eta}^2\big)
\frac{\bm\xi^\dagger\GG\hh_{\ell}^{m}\bm\xi}{\norm{\bm\xi}^2}
=
\frac{1}{\norm{\bm\xi}^{2}}
(1+2\norm{\bm\eta}^2)
\bm\xi^\dagger\GG\hh_{\ell}^{m}\bm\xi.
\]
Comparing the above expression to
\eqref{re-big} with $\bm\eta=0$,
we conclude that
the linearization at a bi-frequency mode
in the invariant subspace
corresponding to spherical harmonics
of degree $\ell$ and orders $\pm m$
is given by the same expression \eqref{lo}
as for one-frequency modes,
but with $(1+2\norm{\bm\eta}^2) m$ in place of $m$,
effectively corresponding to larger values of $m$.
So, if a one-frequency mode
is linearly stable
(with respect to perturbations
in invariant subspaces corresponding to \emph{all}
spherical harmonics),
then a corresponding bi-frequency mode
is also expected to be linearly stable,
at least for $\norm{\bm\eta}$ small enough.
\begin{comment}
Thus, larger values of $\norm{\bm\eta}$
effectively remove harmonics of lower orders $m\ne 0$.
\end{comment}
%Overall, we conclude that
%the linear stability of bi-frequency modes
%is not worse than that of one-frequency modes.

%we have $\jj\p_t R=\mathscr{L}_{\ell,m,\alpha} R$,
%with $\alpha=1+2\sinh\gamma$
%and with $\mathscr{L}_{\ell,m,\alpha}$ given by
\begin{comment}
\[
\mathscr{L}_{\ell,m,\alpha}=
\begin{bmatrix}
L_0&0
\\[2ex]
0&-L_0
\end{bmatrix}
+
\begin{bmatrix}
V&-\alpha m V&V&\alpha m V
\\
0&0&0&0
\\
-V&\alpha m V&-V&-\alpha m V
\\
0&0&0&0
\end{bmatrix}.
\]
\end{comment}

%\section{Numerical results}
\section{Numerical results}
\label{sect-numerics}

We present the spectra
%We computed the real and imaginary part of eigenvalues
of the linearization at a (one-frequency) solitary wave
in invariant spaces $\scrX_{\ell,m}$
for $\abs{m}\leq\ell\le 3$,
given by $\mathbf{A}_{\ell,m}$ from \eqref{lo}.
For simplicity, the mass of the spinor field
is taken $M_s=1$.
%% Below, $m$ denotes the magnetic quantum number
%% (the order of the spherical harmonic),
%% while the mass of spinor field is taken $M_s=1$.
Computation of the spectrum is
similar to \cite{PhysRevLett.116.214101},
but with a differentiation matrix
based on rational Chebysh\"{e}v polynomials in $N=1200$ grid nodes.
We only consider
solitary waves with $\omega\in(0.1,1)$ since as $\omega\to 0$
the numerical accuracy deteriorates
due to the amplitude of solitary waves going to infinity.
%Similarly to what we did for the 2D case, instabilities for $\Im(\lambda)>1-\omega$ are considered to be spurious.
%%The general features of the spectra is that
%the spectrum is symmetric with respect to the real and imaginary axes;
The spectrum of $\mathbf{A}_{\ell,m}$
is symmetric with respect to the real and imaginary axes;
the essential spectrum consists of
$\lambda\in\jj\R$,
$\abs{\lambda}\ge 1-\abs{\omega}$.
\ac{\sout{
  If $\lambda$ is an eigenvalue of a given magnetic quantum number $m$,
  then $\lambda^\dag$ is an eigenvalue for $-m$.
  The union of the spectra for $\pm m$
}}
The spectral (linear) instability is due to eigenvalues
with $\Re\lambda>0$.
% (i.e. the spectrum for $m$ is the complex conjugate of $-m$).

\begin{figure}[ht]
\begin{tabular}{cc}
%\includegraphics[width=4.4cm,height=5cm]{f1.pdf}&
%\includegraphics[width=4.0cm,height=5cm]{f2.pdf}\\
%\quad
%\includegraphics[width=6.6cm,height=7.5cm]{f1.pdf}&
\includegraphics[width=7.2cm,height=8cm]{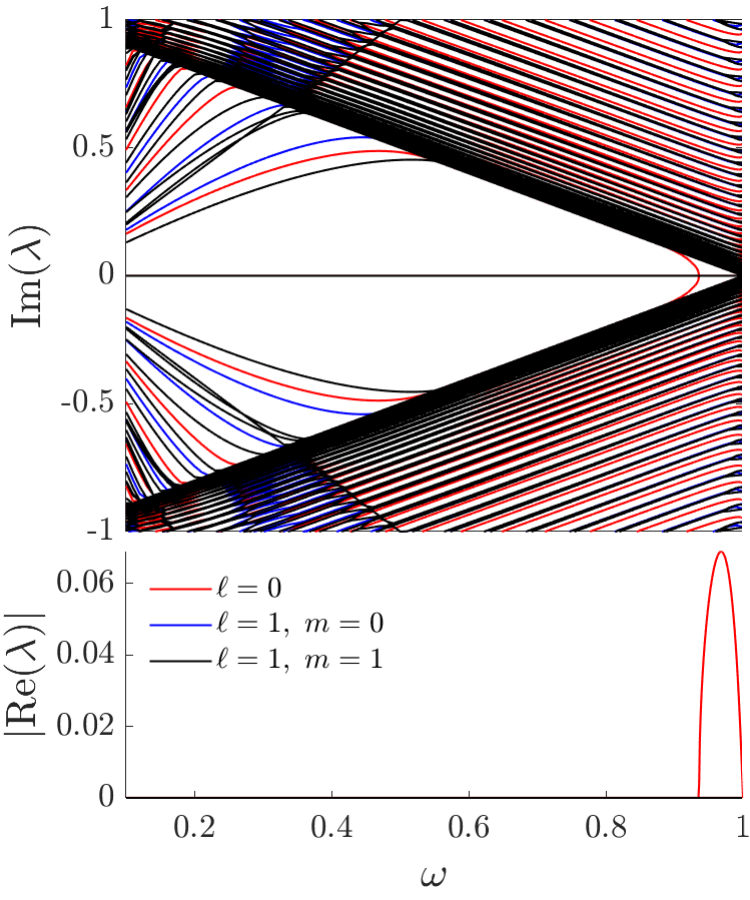}&
\quad
\includegraphics[width=6.6cm,height=8cm]{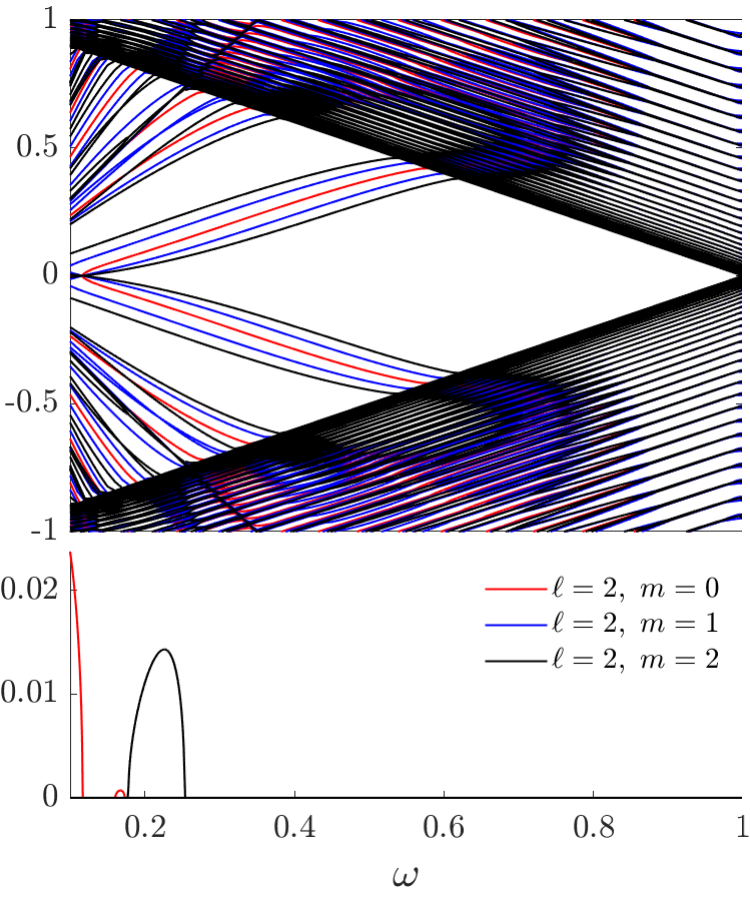}\\
\end{tabular}%
\caption{
Imaginary (top) and real (bottom) parts of the spectrum for
$\ell=0,\,1$ (left) and $\ell=2$ (right)
as functions of $\omega\in(0.1,1)$.
%Note that eigenvalues $0$ and $\pm 2\omega \jj$
%correspond to $\ell=1$, $\abs{m}=1$.
%due to their eigenfunctions being
%highly oscillatory.
}
\label{fig:l012}
\end{figure}

Fig.~\ref{fig:l012}
(left)
shows the spectrum for $\ell=0$ and $\ell=1$.
(Eigenvalue $\lambda=0$ in these cases
corresponds to eigenvectors
$\jj\upvarphi$ and 
$\p_{x_1}\upvarphi$,
$\p_{x_2}\upvarphi$,
$\p_{x_3}\upvarphi$
\cite{berkolaiko2015vakhitov}.)
For $\ell=0$,
the instability region is $\omega\in(0.936,1)$,
due to presence of a pair of real eigenvalues of opposite sign;
these eigenvalues disappear via the pitchfork bifurcation
when $\omega_0\approx 0.936$
\begin{comment}
(which corresponds to solitary waves with minimal charge
\cite{PhysRevD.1.2766}).
\end{comment}
and there are no $\Re\lambda\ne 0$ eigenvalues for $\omega<\omega_0$
\cite{PhysRevLett.116.214101}.
%% A pair of purely imaginary eigenvalues
%% very close to zero
%% for $\omega\lesssim 1$
%% seems to be an artifact:
%% the corresponding eigenfunctions are highly oscillatory.
For $\ell=1$,
there are no $\Re\lambda\ne 0$ eigenvalues;
eigenvalues
$\lambda=\pm2\omega\jj$
stemming from the $\mathbf{SU}(1,1)$ symmetry
\cite{boussaid2018spectral}
correspond to $\abs{m}=\ell=1$.

For $\ell=2$
(right panel of Fig.~\ref{fig:l012}),
for $m=0$,
we found an interval of instability, $\omega\in(0.16,0.174)$,
with a quadruplet of $\Re\lambda\ne 0$ eigenvalues:
this quadruplet appears and disappears
at the endpoints of the interval via the Hamiltonian Hopf (HH) bifurcations,
from the collisions of two pairs of purely imaginary eigenvalues.
(Although the imaginary eigenvalues colliding when $\omega\approx 0.174$
come from the same
threshold, not in line
with the Sturm--Liouville theory expectations,
the form of the eigenfunctions
suggests that this bifurcation is genuine, not a numerical artifact.)
Next onset of instability for $|m|=0$
is from the pitchfork bifurcation at $\omega_p\approx 0.117$.
For $|m|=1$,
there is no instability;
for $|m|=2$,
the instability interval is $\omega\in(0.177,0.254)$,
with the HH bifurcations
at its endpoints.
%possibly also a numerical artifact
%(the same reasons as above).
%
\begin{figure}[ht]
\begin{tabular}{cc}
%\includegraphics[width=4.4cm,height=5cm]{f3.pdf}&
%\includegraphics[width=4.0cm,height=5cm]{f4.pdf}\\
%\quad
\includegraphics[width=7.2cm,height=8cm]{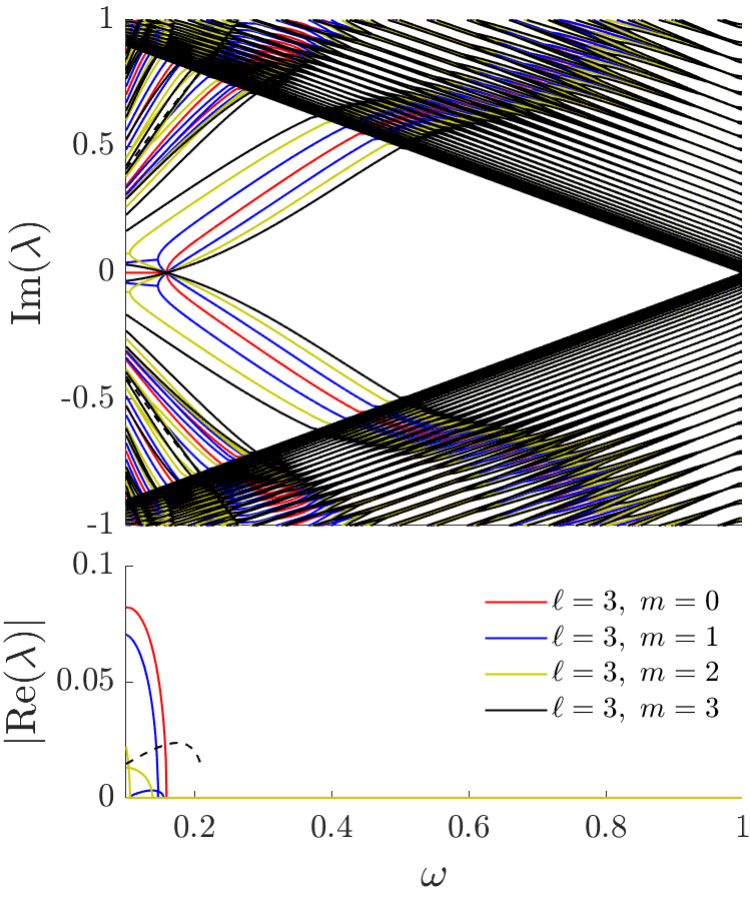}&
\quad
\includegraphics[width=6.6cm,height=8cm]{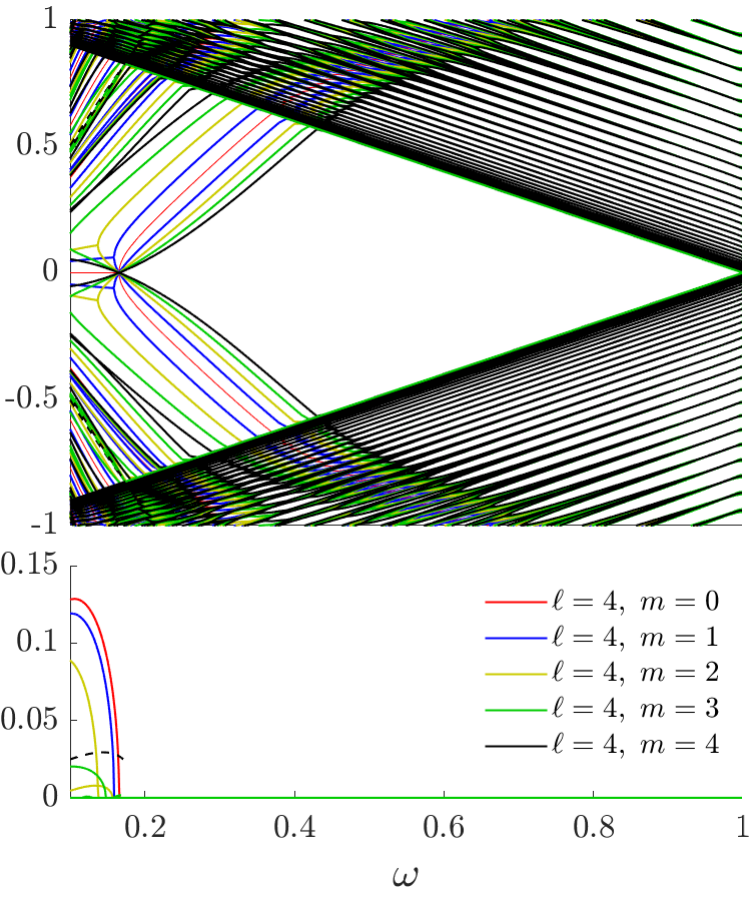}\\
\end{tabular}%
\caption{
%Real and imaginary parts of the spectrum
%as functions of $\omega$
Spectrum for $\ell=3$
(left) and $\ell=4$ (right).
Dashed black lines refer to
quadruplets of eigenvalues
with $\Re\lambda\ne 0$
bifurcating from the thresholds $\pm\jj(1-\omega)$
(possibly a numerical artifact).
}
\label{fig:l34}
\end{figure}

For $\ell=3$
(Fig.~\ref{fig:l34}, left),
%%there are instabilities for every value of $m$.
for $m=0$,
$\Re\lambda>0$ eigenvalue
is born from
the pitchfork bifurcation at $\omega_p\approx 0.159$.
For $|m|=1$, quadruplets of eigenvalues appear
when $\omega$ drops below
$\omega\approx 0.155$
and then below $\omega\approx 0.147$
(the first one disappears at $\omega\approx 0.105$);
for $|m|=2$,
quadruplets appear
at
$\omega\approx 0.139$
and at $\omega\approx 0.106$
(all via HH bifurcations).
For $|m|=3$, there is a quadruplet of $\Re\lambda\ne 0$ eigenvalues
bifurcating from the thresholds $\pm\jj(1-\omega)$
at $\omega\approx 0.2$,
which is
possibly a numerical artifact
since the corresponding eigenfunctions
do not seem to have a continuous limit.

\begin{figure}[ht]
%\begin{tabular}{cc}
%\includegraphics[width=5cm,height=2.4cm]{fig5.pdf}
\begin{center}
\includegraphics[width=8cm,height=6cm]{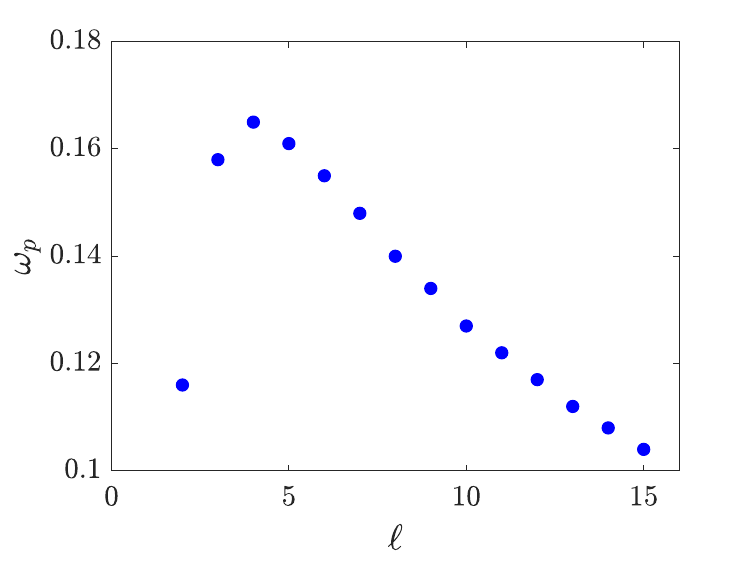}
\end{center}
%\end{tabular}%
\caption{Value $\omega_p$ of the pitchfork bifurcation
for perturbations
from $\scrX_{\ell,m}$.}
\label{fig:5}
\end{figure}

For $\ell=4$ (Fig.~\ref{fig:l34}, right),
for $m=0$,
unstable eigenvalue appears below
pitchfork bifurcation
at $\omega_p\approx 0.166$.
For $|m|=1$, a quadruplet is born at $\omega\approx 0.159$;
for $|m|=2$, another one appears at
$\omega\approx 0.157$ (all via HH bifurcations).
For $|m|=3$, a quadruplet of $\Re\lambda\ne 0$
eigenvalues bifurcating from
the thresholds $\pm\jj(1-\omega)$
when $\omega\approx 0.158$
again seems to be a numerical artifact.
\ac{??}
%there are also overlapping intervals of instability
%$\omega\in(0.116,0.13)$ and $\omega<0.148$.
More quadruplets are born
via HH bifurcations
at $\omega\approx 0.148$
and $\omega\approx 0.13$ (the second disappears at $\omega\approx 0.116$).
%eigenvalues from Hamiltonian Hopf bifurcations
%for $\omega<0.148$.
For $|m|=4$, a quadruplet of $\Re\lambda\ne 0$ eigenvalues
bifurcates from the thresholds
when $\omega\approx 0.17$.
%\ac{????}
%
%% \ac{
%% It can be observed that for $\ell\geq2$ and $m=0$, solitons are exponentially unstable below a critical frequency $\omega_p$. 
%% For $\omega<\omega_p$, there is always a unique pair of 0 eigenvalues. Could it be possible to find for what values of $\omega$ can a 0 eigenvalue pair exist?
%% }
% Another observation is that $\omega_p$ increases with $\ell$.
%What one can wonder is how $\omega_p(\ell)$ grows with $\ell$.
%Does it grow asymptotically or, on the contrary,
% $\omega_p\rightarrow1$ when $\ell\rightarrow\infty$?
%
%\end{comment}
\\
\indent
While the numerics show that
larger $\abs{m}$ lead to smaller intervals
of instability (in agreement with (2+1)D case in \cite{PhysRevLett.116.214101}),
the increase of $\ell$
seems to lead to the growth of the instability interval
$(0,\omega_p(\ell))$.
On Fig.~\ref{fig:5},
one can see that this tendency does not persist:
the maximum value of
$\omega_p\approx 0.166$ occurs for $\ell=4$; for larger $\ell$,
the instability region
$(0,\omega_p(\ell))$ is shrinking.
Let us mention that there is an onset of instability
for $\ell=1$
below the critical value $\omega_p\approx 0.07$,
which is not presented on
Figs.~\ref{fig:l012} and~\ref{fig:5}
since the numerics become less accurate
for small $\omega$.
Thus, the numerics suggest that the spectral stability region
for both one-frequency and bi-frequency modes
is
$\omega\in(0.254 M_s,0.936 M_s)$.

\section{Conclusion}

%We showed that localized modes
%with bi-frequency $e^{\pm\jj\omega t}$ components
%are generically present in
%models of classical spinor fields
%either interacting with scalar fields
%or with self-interaction of scalar type.
%(such as DKG system and the Soler model).
We showed the linear stability of some of
one- and bi-frequency modes
in the (cubic) NLD via the radial reduction.
%We showed that some of bi-frequency modes are also linearly stable
%(and only they can be dynamically stable).
%which indicates
%that the dynamic stability also follows due to
%dispersive properties.
%and showed that such modes are (linearly) stable.
\begin{comment}
Let us point out that until now, stability results
have only been obtained for one-frequency modes.
We stress that in the framework of classical spinor fields
with scalar mechanisms of self-interaction
it is only bi-frequency modes
-- not the one-frequency ones --
that can possess dynamic stability properties
due to dispersion of perturbations.
\end{comment}
%Th approach allows us
%to perform numerical computation of
%computing the spectrum
%of the linearization operator.
%The linear stability approach is based on decomposition
%of a perturbation into spherical harmonics,
%which allows us to perform a radial reduction
%of the linearized equation,
%enabling the numerical computation of the spectrum
% of the linearization operator.
%via via the Evans function approach
%or alternatively 
%via finite difference method.
We presented
the numerical results based on the 
finite difference method,
obtaining a relatively large stability region
$\omega\in(0.254 M_s,0.936 M_s)$,
with $M_s$ the mass of the spinor,
for both one- and bi-frequency modes
in the cubic NLD in (3+1)D.
The perturbation theory implies
that there is a nearby stability region
for one- and bi-frequency
localized modes of the Dirac--Klein--Gordon system
in the case when the boson mass $M_B$
and the coupling constant $g$
%in the DKG system
are large.
%which approaches the above interval
%in the limit of heavy bosons
%and strong interaction.
%and small second term in \eqref{bf-1}.
%Since the charge--current density of bi-frequency spinor modes
%is time-dependent,
Since bi-frequency spinor modes
are only available
in the absence of any interactions
but the Yukawa coupling,
%we suggest that stable bi-frequency modes
we suggest that these modes
can model neutral spinor particles
from the DM sector
which interact with the visible matter
via the ``Higgs portal''.
%(that is, interacting the Higgs field only).

%\begin{floatingfigure}{4.5cm}
%\end{floatingfigure}

\acknowledgments

This research
was supported by a grant from the Simons Foundation (A.C., 851052).
This research has been funded in whole or in part by the French National Research Agency (ANR) as part of the QuBiCCS project ``ANR-24-CE40-3008-01'' (N.B.).
J.C.-M. acknowledges support from the EU (FEDER program 2014--2020) through
MCIN/AEI/10.13039/501100011033 (under the projects PID2020-112620GB-I00 and PID2022-143120OB-I00).

\medskip

\noindent
For the purpose of its open access publication, the author/rights holder applies a CC-BY open access license to any article/manuscript accepted for publication (AAM) resulting from that submission.

%% This is the most common positions for acknowledgments. A macro is
%% available to maintain the same layout and spelling of the heading.

%% \paragraph{Note added.} This is also a good position for notes added
%% after the paper has been written.

\medskip

\noindent{\bf Data Availability Statement.\ }
This article has no associated data or the data will not
be deposited.

\medskip

\noindent{\bf Code Availability Statement.\ }
This article has no associated code or the code will not
be deposited.

\medskip

\noindent{\bf Open Access.\ }
This article is distributed under the terms of the
Creative Commons Attribution License
\href{https://creativecommons.org/licenses/by/4.0/}{(CC-BY4.0)},
 which permits any use, distribution and reproduction in any
medium, provided the original author(s) and source are credited.

\bibliography{bibcomech}
\bibliographystyle{JHEP}
%\bibliographystyle{sima-doi}
%\bibliography{bibcomech}

\end{document}